\documentclass{amsart}
\usepackage{amsaddr}
\usepackage{float}
\usepackage{eucal}
\usepackage[colorlinks,linkcolor=blue,citecolor=blue,urlcolor=blue,linktocpage=true,pdfencoding=auto,unicode,psdextra]{hyperref}
\usepackage[normalem]{ulem}
\usepackage{mathtools}
\usepackage{mathrsfs}
\usepackage{amsmath}
\numberwithin{equation}{section}
\usepackage{lmodern} 
\usepackage{amssymb}
\usepackage{graphicx,scalerel}
\usepackage{xcolor}

\usepackage{calc}
\newsavebox\CBox
\newcommand\hcancel[2][0.5pt]{%
	\ifmmode\sbox\CBox{$#2$}\else\sbox\CBox{#2}\fi%
	\makebox[0pt][l]{\usebox\CBox}%
	\rule[0.5\ht\CBox-#1/2]{\wd\CBox}{#1}}

\def\classificationname{Mathematics Subject Classification}

\newenvironment{classification}[1][2020]{\vspace{.5cm}\small
	\noindent
	\hangindent=\parindent
	\hangafter=1
	{\scshape\classificationname\ (#1).}\ }%
{\par}%

\def\keywordname{Keywords}

\newenvironment{mathkeywords}{\vspace{.3cm}\small
	\noindent
	\hangindent=\parindent
	\hangafter=1
	{\scshape\keywordname.}}%
{\par}%

\def\scientificchapter{Scientific Chapter}

{\par}%

\theoremstyle{plain}
\newtheorem{theorem}{Theorem}[section]

\newtheorem{proposition}[theorem]{Proposition}

\theoremstyle{definition}
\newtheorem{definition}[theorem]{Definition}

\theoremstyle{remark}
\newtheorem{remark} [theorem]{Remark}

\begin{document}
	\title[Weak Lipschitz structures and topological structures]{Weak Lipschitz structures and their connections with the topological structures}
	\author[T. Valent]{Tullio Valent}
	\address{Dipartimento di Matematica ``Tullio Levi Civita'', Universit\`{a} di Padova, Via Trieste 63, 35121 Padova, Italy}
	\email{tullio.valent@unipd.it}

\begin{abstract}
Two approaches to Lipschitz structures for any set $X$ are presented, studied and compared. The first approach is similar to the one proposed in \cite{fras70}, where Lipschitz structures are defined as families of pseudo-metrics satisfying suitable conditions. The other one, here introduced, is expressed by using weak pseudo-metrics, which (unlike the pseudo-metrics) do not necessarily vanish on the whole of the diagonal of $X \times X$; in this case we will talk about weak Lipschitz structures. Since all topological structures are defined by a a family of \emph{weak pseudo-metrics} (as we will show in Section \ref{sect4}) we can find some connections between topological structures and weak Lipschitz structures, and a link between continuous maps and weak Lipschitz maps.

A central part of this paper is devoted to the \emph{weak Lipschitz uniformity} defined by a weak Lipschitz structure, which is introduced in Section \ref{sect8}. A notion of uniform continuity with respect to weak Lipschitz uniformities is proposed and studied. In particular, we prove that the weak Lipschitz maps acting between two weak Lipschitz spaces $(X,  \mathscr{L}_X)$ and $(Y,  \mathscr{L}_Y)$ are uniformly continuous with respect to the weak Lipschitz uniformities defined by $\mathscr{L}_X$ and $\mathscr{L}_Y$.
\end{abstract}

\maketitle

\begin{classification}
	Primary: 54A05; Secondary: 54C08.
\end{classification}

\begin{mathkeywords}
	{General structure theory, Lipschitz structures and topologies}
\end{mathkeywords}

\section{Introduction.}\label{intro}
With the aim of studying Lipschitz structures for sets endowed with some other structures (in particular for topological vector spaces), in this paper we consider Lipschitz and weak Lipschitz structures for any set.

We will present and compare two types of Lipschitz structures for a set $X$: the Lipschitz structures originally presented in \cite{fras70}, and the \emph{weak Lipschitz structures} here introduced, which are a suitable family of \emph{weak pseudo-metrics}. By a weak pseudo metric on a set $X$ we mean a symmetric map  $d: X \times X \mapsto \mathbb{R}^+$ such that $d(x_1,x_2) \leq d(x_1,x) + d(x,x_2)$ for all $x_1,x_2,x \in X$, which vanishes in at least a point of the diagonal $\Delta$ of $X \times X$ (but not necessarily on the whole of $\Delta$).

In Section \ref{sect4} we exhibit a link between topological structures and weak Lipschitz structures. Section \ref{sect5} is devoted to the study of products of Lipschitz spaces and of weak Lipschitz spaces, while in Section \ref{sect6} we define the locally Lipschitz maps and the locally weak Lipschitz maps.

The subsequent Sections \ref{sect7} and \ref{sect8} contain the main results of this paper, which concern Lipschitz and weak Lipschitz uniformities. In them we prove, among other results, the uniform continuity of the weak Lipschitz maps with respect to the weak Lipschitz uniformities.


\section{Lipschitz structures and weak Lipschitz structures. Lipschitz maps and weak Lipschitz maps.}\label{sect2}
Let $X$ be any set. Any symmetric map $d : X \times X \mapsto \mathbb{R}^+$ such that $d(x_1,x_2) \leq d(x_1,x) + d(x,x_2) \quad \forall x,x_1,x_2 \in X$ and   vanishes in at least a point of the diagonal $\Delta$ of $X \times X$ will be called a \emph{weak pseudo-metric} on $X$. We emphasize the fact
that $d$ is not requested to vanish on the whole of $\Delta$. The family of all pseudo-metrics on $X$ will by denoted by $\mathscr{P}(X)$, and the family of all weak pseudo-metrics on $X$ will by denoted by $\mathscr{P}_w(X)$.

\begin{definition}
	By a \emph{Lipschitz structure} (resp. \emph{weak Lipschitz structure}) for $X$   we mean a non-empty family $\mathscr{L}$ of pseudo-metrics (resp. weak pseudo-metrics) on $X$  satisfying the following conditions:
	
	\begin{itemize}
		\item[$(L_1)$] $d \leq d_1, d_1 \in \mathscr{L} \Rightarrow d \in  \mathscr{L}$;
		\item[$(L_2)$] $d \in \mathscr{L} \Rightarrow \alpha d \in  \mathscr{L}$ for every real number $\alpha >0$;
		\item[$(L_3)$] $d_1,d_2 \in \mathscr{L} \Rightarrow d_1 \vee d_2 \in \mathscr{L}$.
	\end{itemize}

\end{definition}

Of course, a Lipschitz structure is a particular weak Lipschitz structure.
It is easy to prove

\begin{remark}\label{rem1}
	$\left((L_1),(L_2),(L_3)\right)$ is equivalent to 	$\left((L_1),(L_4)\right)$, where 
	\begin{itemize}
		\item[$(L_4)$] $d_1,d_2 \in \mathscr{L} \Rightarrow d_1 + d_2 \in \mathscr{L}$.
	\end{itemize}
	
\end{remark}

The pair $(X,\mathscr{L})$ will be called a \emph{Lipschitz space} when $\mathscr{L}$ is a Lipschitz structure and a \emph{weak Lipschitz space} if $\mathscr{L}$ is a weak Lipschitz structure.

\begin{definition}
	\emph{A base of a (weak) Lipschitz structure $\mathscr{L}$} is a subset $\mathscr{B}$ of $\mathscr{L}$ such that for every $d \in \mathscr{L}$ there are $b \in \mathscr{B}$ and $\alpha >0$ such that
	$$d \leq \alpha b.$$
\end{definition}

If $\mathscr{B}$ is a base of $\mathscr{L}$ then 
$$\mathscr{L} = \{d: d \leq \alpha b \quad \text{for some} \quad b \in \mathscr{B} \quad \text{and}\quad \alpha >0\}.$$

\begin{definition}
	A set $\mathscr{B}$ of pseudo-metrics (resp. of weak pseudo-metrics) on $X$ is \emph{a base for a Lipschitz structure (resp. for a weak Lipschitz structure) for $X$}  if there is a Lipschitz structure (resp. a weak Lipschitz structure) for $X$ of which $\mathscr{B}$ is a base. This is true if and only if the set
	$$\{d: d \leq \alpha b \quad \text{with} \quad b \in \mathscr{B} \quad \text{and}\quad \alpha >0\}$$
	is a Lipschitz structure (resp. weak Lipschitz structure) for $X$.
\end{definition}
It follows that $\mathscr{B}$ is a base for a Lipschitz structure, or for a weak Lipschitz structure,  for $X$ if and only if  for every $b_1,b_2 \in \mathscr{B}$ there are $b \in \mathscr{B}$ and $\alpha >0$ such that $b_1 \vee b_2 \leq \alpha b$. Consequently, given a non-empty family $\mathscr{P}$ of pseudo-metrics (resp. $\mathscr{P}_w$ of weak pseudo-metrics) on $X$, the family $\mathscr{B}$ of the suprema of all finite subsets of $\mathscr{P}$ (resp. $\mathscr{P}_w$) is a base for a Lipschitz structure $\mathscr{L}(\mathscr{P})$ (resp. for a weak Lipschitz structure $\mathscr{L}(\mathscr{P}_w))$ for $X$. Of course, $\mathscr{L}(\mathscr{P})$ is the smallest Lipschitz structure for $X$ containing $\mathscr{P}$, and it will be called the \emph{weak Lipschitz structure for X generated by $\mathscr{P}$}, while $\mathscr{L}(\mathscr{P}_w)$ is the smallest weak Lipschitz structure for $X$ containing $\mathscr{P}_w$, and it will be called the \emph{weak Lipschitz structure for X generated by $\mathscr{P}_w$}. 

We have
$$\mathscr{L}(\mathscr{P})=\{d \in \mathscr{P}(X): d \leq \alpha (d_1\vee  \dots \vee d_n), d_1,\dots,d_n \in \mathscr{P}, n \geq 1, \alpha >0 \}$$  
and
$$\mathscr{L}(\mathscr{P}_w)=\{d \in \mathscr{P}_w(X): d \leq \alpha (d_1\vee  \dots \vee d_n), d_1,\dots,d_n \in \mathscr{P}_w, n \geq 1, \alpha >0 \}.$$ 

\begin{definition}
	Let $(X,\mathscr{L}_X)$ and $(Y,\mathscr{L}_Y)$ be Lipschitz spaces (resp. weak Lipschitz spaces). A map $f: X \mapsto Y$ is called a \emph{Lipschitz map}  (resp. \emph{weak Lipschitz map}) if for every $d_Y \in \mathscr{L}_Y$ there is $d_X \in \mathscr{L}_X$ such that
	
	$$d_Y(f(x_1),f(x_2))\leq d_X(x_1,x_2) \quad \forall x_1,x_2 \in X,$$
	i.e.,
	$$d_Y \circ \overset{(2)}{f} \leq d_X.$$
	
\end{definition}	

\begin{remark}\label{rem:1}
	Let $\mathscr{L}_Y$ be generated by a family  $\mathscr{P}_Y$ of pseudo-metrics (resp. of weak pseudo-metrics) on $Y$. A map $f:X \mapsto Y$ is a Lipschitz map (resp. a weak Lipschitz map) if for every $d \in \mathscr{P}_Y$ there is $d_X \in \mathscr{L}_X$ such that $d \circ \overset{(2)}{f} \leq d_X.$
\end{remark}

\begin{proof}
	Suppose that  for every $d \in \mathscr{P}_Y$ there is $d_X \in \mathscr{L}_X$ such that $d \circ \overset{(2)}{f} \leq d_X.$ We must prove that, consequently, for every $d_Y \in \mathscr{L}_Y$ there is $\delta_X \in \mathscr{L}_X$ such that  $d_Y \circ \overset{(2)}{f} \leq \delta_X$. Indeed, given $d_Y \in \mathscr{L}_Y$, since $\mathscr{P}_Y$ generates $\mathscr{L}_Y$ there are a finite subset $\{d_1,\dots,d_n\}$ of $\mathscr{P}_Y$ and a number $\alpha >0$ such that $d_Y \leq \alpha (d_1\vee  \dots \vee d_n)$. Let $\delta_1,\dots,\delta_n$ be elements of $\mathscr{L}_X$ such that $d_1 \circ \overset{(2)}{f} \leq \delta_1,\dots, d_n \circ \overset{(2)}{f} \leq \delta_n$. It follows
	$$d_Y \circ \overset{(2)}{f} \leq \left(\alpha(d_1\vee  \dots \vee d_n)\right) \circ \overset{(2)}{f} = 
	\alpha ((d_1 \circ  \overset{(2)}{f}) \vee \dots \vee (d_n \circ  \overset{(2)}{f})) \leq \alpha (d_X^1 \vee \dots \vee d_X^n),$$
	with $d_X^1,\dots,d_X^n$ elements of $\mathscr{P}_Y$. Since \mbox{$\alpha (d_X^1 \vee \dots \vee d_X^n) \in \mathscr{L}_X$} the proof is concluded.
\end{proof}

\section{Lipschitz structures and weak Lipschitz structures for $X$ defined by a family of subsets of $X$.}\label{sect3}

Let $\mathscr{A}$ be a non empty family of subsets of $X$, and let us set
\begin{equation*}
	\begin{aligned}
		\mathscr{L}_1(\mathscr{A}) &\coloneqq \{d \in \mathscr{P}(X): d \,\; \text{is bounded on} \, A \times A \,\; \text{for some} \, A \in \mathscr{A}\},\\
		\mathscr{L}_2(\mathscr{A}) &\coloneqq \{d\in \mathscr{P}(X): d \,\;  \text{is bounded on} \, A \times A \,\; \text{for all} \, A \in \mathscr{A}\}.\\
	\end{aligned}
\end{equation*}

\begin{theorem}\label{thm:1}
	If $\mathscr{A}$ is closed under finite intersections, then $\mathscr{L}_1(\mathscr{A})$ is a  Lipschitz structure.
\end{theorem}
\begin{proof}
	Obviously, if $d_1 \in \mathscr{L}_1(\mathscr{A})$ and $d$ is a  pseudo-metric on $X$ smaller than $d_1$, then $ d \in \mathscr{L}_1(\mathscr{A})$. We now show that 
	$$d_1,d_2 \in \mathscr{L}_1(\mathscr{A}) \Rightarrow d_1 + d_2 \in \mathscr{L}_1(\mathscr{A}).$$
	This is true, because if $A_1,A_2 \in \mathscr{A}$ and 
	$$ \sup d_1(A_1 \times A_1) = a_1 \in \mathbb{R}, \quad \sup d_2(A_2 \times A_2) = a_2 \in \mathbb{R},$$
	then
	$$\sup (d_1+d_2)(A \times A) \leq a_1 + a_2,\quad \text{with} \quad A= A_1 \cap A_2.$$
\end{proof}

\begin{theorem}\label{thm:2}
	For any non-empty family of subsets of $X$, $\mathscr{L}_2(\mathscr{A})$ is a Lipschitz structure.
\end{theorem}
\begin{proof}
	The implication $d \leq d_1 \in \mathscr{L}_2(\mathscr{A}) \Rightarrow d \in \mathscr{L}_2(\mathscr{A})$ is evident. The implication $d_1,d_2 \in \mathscr{L}_2(\mathscr{A}) \Rightarrow d_1 + d_2 \in \mathscr{L}_2(\mathscr{A})$ is easily showed. Indeed, if $d_1,d_2 \in \mathscr{L}_2(\mathscr{A})$ there  is $A \in \mathscr{A}$ such that $\sup d_1(A \times A) \in \mathbb{R}$ and $\sup d_2(A \times A) \in \mathbb{R}$. Consequently $\sup (d_1 + d_2)(A \times A) \in \mathbb{R}$, and this implies $ d_1 + d_2 \in \mathscr{L}_2(\mathscr{A})$.
\end{proof}

\begin{remark}
	Theorem \ref{thm:1} and Theorem \ref{thm:2} hold also for weak Lipschitz structures, with 
	
	\begin{equation*}
		\begin{aligned}
			\mathscr{L}_1(\mathscr{A}) &\coloneqq \{d \in \mathscr{P}_w(X): d \,\; \text{is bounded on} \, A \times A \,\; \text{for some} \, A \in \mathscr{A}\},\\
			\mathscr{L}_2(\mathscr{A}) &\coloneqq \{d\in \mathscr{P}_w(X): d \,\;  \text{is bounded on} \, A \times A \,\; \text{for all} \, A \in \mathscr{A}\}.\\
		\end{aligned}
	\end{equation*}
\end{remark}

\section{Connections between topological structures and weak Lipschitz structures}\label{sect4}
Through this paper, when   $\mathscr{P}$ is a family of weak pseudo-metrics on $X$, we will consider as \emph{the topology defined by $\mathscr{P}$}  the one generated by the family of the subsets $U_{d,\varepsilon}(x)$ of $X$, with $d \in \mathscr{P}, x \in X,$ and $ \varepsilon$ a real number such that $\varepsilon>d(x,x)$, where
\begin{equation*}
	U_{d,\varepsilon}(x) \coloneqq \{\xi\in X : d(\xi,x)<\varepsilon\}.
\end{equation*}	

\begin{theorem}\label{thm:topdef}
	Every topology $\tau$ on $X$ is defined by the family 
	$$\mathscr{P}_{\tau}\coloneqq \{d_A: A \in \tau\},$$
	where $d_A$ is the weak pseudo-metric defined as the  characteristic function of the complement \mbox{$\complement(A \times A)$} of $A \times A$. 
\end{theorem}
\begin{proof}
	It is easy to see that the maps $d_A: X \times X \mapsto \mathbb{R}^+$ are weak pseudo-metrics. 
	The topology defined by $\mathscr{P}_{\tau}$ is the topology generated by the family of the subsets $U_{d_{A},\varepsilon}(x)$ of $X$, where $A \in \tau, x \in X, \varepsilon>d_{A}(x,x)$, and 
	\begin{equation*}
		U_{d_A,\varepsilon}(x) \coloneqq \{\xi\in X : d_A(\xi,x)<\varepsilon\}.
	\end{equation*}
	So, the topology defined by the family $\mathscr{P}_{\tau}$ of weak pseudo-metrics coincides with $\tau$, because $U_{d_A,\varepsilon}(x)=A$ if $x \in A$ and $\varepsilon \leq 1$, $U_{d_A,\varepsilon}(x)=X$ if $x \in A$ and $\varepsilon > 1$, while $ U_{d_A,\varepsilon}(x)=X$ if $x \notin A$ and $\varepsilon > 1 $, and $U_{d_A,\varepsilon}(x)=\emptyset$ when $x \notin A$ and $\varepsilon \leq 1$. 
\end{proof}

\begin{definition}
For every weak Lipschitz structure $\mathscr{L}$ for a set $X$ we will denote by $\widehat{\tau}(\mathscr{L})$ \emph{the topology on $X$ defined by the family  $\mathscr{L}$ of weak pseudo-metrics}.
\end{definition}
\begin{definition}
For every topology $\tau$ on a set $X$ we will denote by  $\widehat{\mathscr{L}}(\tau)$ \emph{the smallest weak Lipschitz structure for X containing the family $\mathscr{P}_{\tau}$ of weak pseudo-metrics} considered in Theorem \ref{thm:topdef}. In other words,  $\widehat{\mathscr{L}}(\tau)$ is \emph{the weak Lipschitz structure for $X$ generated by  $\mathscr{P}_{\tau}$}.
\end{definition}

\begin{theorem}
	Let $f:X \mapsto Y$, and let $\tau_X, \tau_Y$ be topologies on $X$ and $Y$. If $f$ is continuous, then $f$ is a weak Lipschitz map with respect to the weak Lipschitz structures $\widehat{\mathscr{L}}(\tau_X)$ for $X$ and $\widehat{\mathscr{L}}(\tau_Y)$ for $Y$. 
\end{theorem}

\begin{proof}
	Let $f$ be continuous for the topologies $\tau_X$ on $X$ and  $\tau_Y$ on  $Y$. We will show that for every $B \in \tau_Y$ there is $A \in \tau_X$ such that 
	\begin{equation}\label{funceq}
		d_B(f(x_1),f(x_2)) \leq d_A(x_1,x_2) \quad \forall x_1,x_2 \in X,
	\end{equation} 
	so proving (in view of Remark \ref{rem:1}) that $f$ is a weak Lipschitz map with respect to the weak Lipschitz structures $\widehat{\mathscr{L}}(\tau_X)$  and $\widehat{\mathscr{L}}(\tau_Y)$. Set, for every $B \in \tau_Y$, $A \coloneqq f^{\leftarrow}(B)$. $A$ is an element of $\tau_X$ because $f$ is continuous. The inequality \eqref{funceq} is obviously satisfied when $(x_1,x_2) \in \complement(A \times A)$. If $(x_1,x_2) \in A \times A$, then $(f(x_1),f(x_2)) \in  B \times B$, which implies $d_B(f(x_1),f(x_2))=0$. Hence we can conclude that \eqref{funceq} is true for all $x_1,x_2 \in X$.
\end{proof} 

\begin{theorem}
	Let  $f$ be a weak Lipschitz map with respect to the weak Lipschitz structures $\mathscr{L}_X$ for $X$ and $\mathscr{L}_Y$ for $Y$. Then $f$ is continuous for the topologies $\widehat{\tau}(\mathscr{L}_X)$ on $X$ and $\widehat{\tau}(\mathscr{L}_Y)$ on $Y$.
\end{theorem}

\begin{proof}
	We recall that $\widehat{\tau}(\mathscr{L}_X)$ is the topology on $X$ defined by the family $\mathscr{L}_X$ of weak pseudo-metrics, and $\widehat{\tau}(\mathscr{L}_Y)$ is the topology on $Y$ defined by the family $\mathscr{L}_Y$ of weak pseudo-metrics.
	
	We have supposed that for every weak pseudo-metric $d_Y \in \mathscr{L}_Y$ there is a weak pseudo-metric  $d_X \in \mathscr{L}_X$ such that
	$$d_Y(f(x_1),f(x_2)) \leq d_X(x_1,x_2) \quad \forall x_1,x_2 \in X.$$
	It follows that, for every $x \in X$, 
	$$d_Y(f(\xi),f(x)) \leq d_X(\xi,x) \quad \forall \xi \in X.$$
	Thus we can conclude that $f$ is continuous at any $x \in X$ for the topologies $\widehat{\tau}(\mathscr{L}_X)$ on $X$ and $\widehat{\tau}(\mathscr{L}_Y)$ on $Y$.
\end{proof}

\section{Product of Lipschitz spaces and of weak Lipschitz spaces.}\label{sect5}

Let $(X_i,\mathscr{L}_i)$, $i \in I$, be a finite family of Lipschitz spaces. Consider the maps \\ $ d: (\prod_{i \in I} X_i) \times (\prod_{i \in I} X_i) \mapsto \mathbb{R}^+$ defined by setting, for $(x_i)_{i \in I}, (\xi_i)_{i \in I} \in \prod_{i \in I} X_i$,

\begin{equation}\label{eq:distproduct}
	d((x_i)_{i \in I},(\xi_i)_{i \in I}) = \sum_{i \in I} d_i(x_i,\xi_i),
\end{equation}

with $d_i$ any element of $\mathscr{L}_i$.

\begin{proposition}\label{prop:product}
	The maps $d$ defined by \eqref{eq:distproduct} are pseudo-metrics on the product $\prod_{i \in I} X_i$.
	\end{proposition}

\begin{proof}It will be shown the triangle inequality:
\begin{equation*}
	d((x_i)_{i \in I},(\xi_i)_{i \in I}) \leq d((x_i)_{i \in I},(\eta_i)_{i \in I}) + d((\eta_i)_{i \in I},(\xi_i)_{i \in I}),
\end{equation*}
namely 

\begin{equation*}
	\sum_{i \in I} d_i(x_i,\xi_i) \leq \sum_{i \in I} d_i(x_i,\eta_i) + \sum_{i \in I} d_i(\eta_i,\xi_i)
\end{equation*}
for all $(x_i)_{i \in I}, (\xi_i)_{i \in I}, (\eta_i)_{i \in I} \in \prod_{i \in I} X_i$.

This inequality is satisfied because the maps $d_i$ are pseudo-metrics and so, for every $i \in I$, we have
\begin{equation*}
	d_i(x_i,\xi_i) \leq d_i(x_i,\eta_i) + d_i(\eta_i,\xi_i).
\end{equation*}
The proof is concluded, since the maps $d$ obviously have the others properties required to be pseudo-metrics.
\end{proof}

\begin{definition}
	The \emph{product of the Lipschitz spaces} $(X_i,\mathscr{L}_i),  i \in I$, is the set  $\prod_{i \in I} X_i$ endowed with the Lipschitz structure generated by the family of the pseudo-metrics $d$ defined by \eqref{eq:distproduct}, i.e. the smallest Lipschitz structure containing all the pseudo-metrics $d$ defined by \eqref{eq:distproduct}.
	\end{definition}

Let now $(X_i,\mathscr{L}_i^w),  i \in I$ be a finite family of weak Lipschitz spaces, and consider the maps  $ d: (\prod_{i \in I} X_i) \times (\prod_{i \in I} X_i) \mapsto \mathbb{R}^+$ defined by \eqref{eq:distproduct} with $d_i \in 
\mathscr{L}_i^w$. An argument similar to the one used in the proof of Proposition \ref{prop:product} shows that the maps $d$ are weak pseudo-metrics on $\prod_{i \in I} X_i$. This leads us to give the following

\begin{definition}
	The \emph{product of the weak Lipschitz spaces} $(X_i,\mathscr{L}_i^w),  i \in I$, is the set  $\prod_{i \in I} X_i$ endowed with the weak Lipschitz structure generated by the family of the weak pseudo-metrics $d$ defined by \eqref{eq:distproduct} with  $d_i \in 
	\mathscr{L}_i^w$.
	\end{definition}
	
\section{Locally Lipschitz maps. Locally scalar Lipschitz maps.}\label{sect6}	
Let $(X,\mathscr{L}_X)$ and $(Y,\mathscr{L}_Y)$	be Lipschitz spaces, or weak Lipschitz spaces. 

\begin{definition}
	A map $f: X \mapsto Y$ will be called a \emph{locally Lipschitz map} (resp. \emph{a locally weak Lipschitz map}) if every $x \in X$ has a neighborhood $U_x$ for the topology associated to $\mathscr{L}_X$ such that for every $d_Y$ belonging to the Lipschitz (resp. weak Lipschitz) structure $\mathscr{L}_Y$ there is an element $d_X$ of the Lipschitz (resp. weak Lipschitz) structure $\mathscr{L}_X$ such that
	
	  	\begin{equation*}\label{funceq2}
	  	d_Y(f(x_1),f(x_2)) \leq d_X(x_1,x_2) \quad \forall x_1,x_2 \in U_x.
	  \end{equation*} 
  \end{definition}

Recalling how a topology may be associated to a Lipschitz or to a weak Lipschitz structure, it is easy to check the following remark.

\begin{remark}
	A map  $f: X \mapsto Y$ is a locally Lipschitz map if for each $d_Y \in \mathscr{L}_Y$ and each $x \in X$ there are  $d_X \in \mathscr{L}_X$ and a number $r_x > d_X(x,x)$ such that $d_X(x,\xi) < r_x$ $\Rightarrow$ $d_Y(f(x),f(\xi)) < d_X(x,\xi)$ for any $\xi \in X$.
	\end{remark}

\begin{definition}
	A map $f: X \mapsto Y$ will be called a \emph{(locally) scalar Lipschitz map} [resp. \emph{a (locally) scalar weak Lipschitz map}] if for every (locally) Lipschitz map [resp. weak Lipschitz map] $\varphi : Y \mapsto \mathbb{R}$ the map $\varphi \circ f : X \mapsto \mathbb{R}$ is a (locally) Lipschitz map [resp. a (locally) weak Lipschitz map].
	\end{definition}

\section{Lipschitz uniformities. Uniform continuity of Lipschitz maps.}\label{sect7}
Let us recall (see, as example, \cite{kell55}) that every family $\mathscr{P}$  of  pseudo-metrics for a set $X$ defines a \emph{uniformity} $\mathscr{U}(\mathscr{P})$ for $X$. 
$\mathscr{U}(\mathscr{P})$ is the filter on $X \times X$ 
generated by the family of the sets $\{(x_1,x_2) \in X \times X : \, d(x_1,x_2) < \varepsilon\}$, with $d \in \mathscr{P}$ and $\varepsilon$ a positive real number.

A function $f$ from a uniform space $(X,\mathscr{U})$ into an uniform space $(Y,\mathscr{V})$ is said to be \emph{uniformly continuous} if for each $V \in \mathscr{V}$ the set  $\{(x,y) \in X \times Y : (f(x),f(y)) \in V\}$ is an element of $\mathscr{U}$.

The \emph{product uniformity} of $\mathscr{U}$ and $\mathscr{V}$ is the smallest uniformity for $X\times Y$ such that the projections into $(X,\mathscr{U})$ and into $(Y,\mathscr{V})$ are uniformly continuous. It is known that, if $(X,\mathscr{U})$ is an uniform space, a pseudo-metric $d : X \times X \mapsto \mathbb{R}^+$ is uniformly continuous with respect to the product uniformity on $X \times X$ and the usual uniformity on $\mathbb{R}$ if and only if the set $\{(x_1,x_2) \in X \times X : \, d(x_1,x_2) < \varepsilon\}$ is an element of $\mathscr{U}$ for each positive real number $\varepsilon$.

If $\mathscr{U}(\mathscr{P})$ is the uniformity for $X$ defined by a family $\mathscr{P}$ of pseudo-metrics, then these pseudo-metrics are uniformly continuous. An important theorem concerning the uniformity theory asserts that \emph{each uniformity for $X$ is defined by the family of all pseudo-metrics for $X$ which are uniformly continuous}.

It follows that the largest family of pseudo-metrics which defines an uniformity $\mathscr{U}$ for $X$ is the family of all pseudo-metrics which are uniformly continuous with respect to the product uniformity of $\mathscr{U}$ by itself. Let us now apply these facts to the Lipschitz structures. 

\begin{definition}
	If $\mathscr{L}$ is a Lipschitz structure, the uniformity $\mathscr{U}(\mathscr{L})$ defined by the family $\mathscr{L}$ of pseudo-metrics for $X$ will be said a \emph{Lipschitz uniformity}.
	\end{definition}

It easily follows (keeping in mind the definition of Lipschitz map) the following

\begin{theorem}
	Each Lipschitz map is uniformly continuous with respect to the Lipschitz uniformities.
	\end{theorem}

\begin{proposition}
	The elements of $\mathscr{L}$ are (pseudo-metrics) uniformly continuous on $X \times X$ with respect to the product of the Lipschitz uniformity $\mathscr{U}(\mathscr{L})$ by itself on $X \times X$ and the usual uniformity on $\mathbb{R}$.
	\end{proposition}

\begin{remark}
	The largest family of pseudo-metrics defining the Lipschitz uniformity  $\mathscr{U}(\mathscr{L})$ is, in general, wider that $\mathscr{L}$.
	\end{remark}

\section{Weak Lipschitz uniformity defined by a weak Lipschitz structure. Uniform continuity of weak Lipschitz maps.}\label{sect8}
\begin{proposition}
	Let $\mathscr{L}$ be a weak Lipschitz structure for $X$. The family of the subsets $U_{d,\varepsilon}$ of $X \times X$ defined by
	\begin{equation*}
		U_{d,\varepsilon} \coloneqq \{(x_1,x_2) \in X \times X: \, d(x_1,x_2)< \varepsilon\},
	\end{equation*}
with $d \in \mathscr{L}$ and $\varepsilon$ a positive real number, is a pre-base for a filter $\mathscr{U}(\mathscr{L})$ on $X \times X$.
\end{proposition}
\begin{proof}
	We must prove that for each $d_1,d_2 \in \mathscr{L}$ and $\varepsilon_1, \varepsilon_2$ positive real numbers the intersection $	U_{d_1,\varepsilon_1} \cap \,	U_{d_2,\varepsilon_2}$ is non empty. This is true because 
	 $	U_{d_1,\varepsilon_1} \cap 	U_{d_2,\varepsilon_2}= \{(x_1,x_2) \in X \times X: \, d_1(x_1,x_2)< \varepsilon_1, 
	 d_2(x_1,x_2)< \varepsilon_2\}$ contains the set $\{(x_1,x_2) \in X \times X: \, (d_1 \vee d_2)(x_1,x_2)< \varepsilon_1 \wedge \varepsilon_2\}$ which is non-empty since the element $d_1 \vee d_2$ of $\mathscr{L}$ vanishes at some point of the diagonal of $X \times X$. 
\end{proof}

\begin{definition}
	The filter $\mathscr{U}(\mathscr{L})$ will be called the \emph{weak Lipschitz uniformity for $X$ defined by $\mathscr{L}$}.
	\end{definition}
There is a (at least formal) analogy between some properties of the pseudo-uniformities introduced and studied in \cite{val23} and those of the weak Lipschitz uniformities here considered. Observe, for example, that they are both defined by a family of (uniformly continuous) weak pseudo-metrics.
\begin{definition}
	Let $(X, \mathscr{L}_X)$ and $(Y, \mathscr{L}_Y)$ be weak Lipschitz spaces. A map $f: X \mapsto Y$ will be called \emph{uniformly continuous with respect to the weak Lipschitz uniformities $\mathscr{U}(\mathscr{L}_X)$ and $\mathscr{U}(\mathscr{L}_Y)$} if for each $U \in \mathscr{U}(\mathscr{L}_Y)$ the set $\{(x_1,x_2) \in X \times X: \, (f(x_1),f(x_2)) \in U \}$ is an element of $\mathscr{U}(\mathscr{L}_X)$.
	\end{definition}

\begin{theorem}
	Let $(X, \mathscr{L}_X)$ and $(Y, \mathscr{L}_Y)$ be weak Lipschitz spaces. A weak Lipschitz map $f: X \mapsto Y$ is uniformly continuous with respect to the weak Lipschitz uniformities defined by $\mathscr{L}_X$ for $X$ and $\mathscr{L}_Y$ for $Y$.
	\end{theorem}
\begin{proof}
	Let $U_Y \in \mathscr{U}(\mathscr{L}_Y)$. Then there are a finite subset $\{d_1,\dots,d_n\}$ of $\mathscr{L}_Y$ and positive real numbers $r_1,\dots,r_n$ such that 
	\begin{equation*}
		y_1,y_2 \in Y, \: d_i(y_1,y_2) < r_i \: \: \forall i=1,\dots,n \Rightarrow (y_1,y_2) \in U_Y.
	\end{equation*}
Since $f$ is a weak Lipschitz map there are $\delta_1,\dots, \delta_n \in \mathscr{L}_X$ such that 
\begin{equation*}
	d_i(f(x_1),f(x_2)) \leq \delta_i(x_1,x_2) \quad \text{for all} \quad x_1,x_2 \in X,
\end{equation*}
which implies $	d_i(f(x_1),f(x_2)) < r_i$ whenever $\delta_i(x_1,x_2) < r_i$.
If we set $$U_X \coloneqq \{(x_1,x_2) \in X \times X : \delta_i(x_1,x_2) < r_i \: \: \forall i=1,\dots,n\}$$ we have $U_X \in \mathscr{U}(\mathscr{L}_X)$, and $(f(x_1),f(x_2)) \in U_Y$ for all $(x_1,x_2) \in U_X$, namely
\begin{equation*}
	U_X \subseteq \{(x_1,x_2) \in X \times X : (f(x_1),f(x_2)) \in U_Y\}.
\end{equation*}
Thus
\begin{equation*}
	\{(x_1,x_2) \in X \times X : (f(x_1),f(x_2)) \in U_Y\} \in \mathscr{U}(\mathscr{L}_X).
\end{equation*}
Therefore we can conclude that $f$ is uniformly continuous with respect to the weak Lipschitz uniformities $\mathscr{U}(\mathscr{L}_X)$ and $\mathscr{U}(\mathscr{L}_Y)$.
\end{proof}

%
%
%
%
%
%


\footnotesize
\begin{thebibliography}{99}
\bibitem{fras70} \textsc{Fraser, Jr.~R. B.}
\textit{Axiom systems for Lipschitz structures,}
Fundamenta Mathematicae, {\bf 66} 1 (1970), 15--24.

%
\bibitem{kell55}\textsc{Kelley~J. L.}, 
{\it General Topology},
D. Van Nostrand Co. Inc, Princeton, N. J., 1955.

\bibitem{val23} \textsc{Valent, T.}
\textit{Pseudo-uniformities,}
Atti Accad. Naz. Lincei Cl. Sci. Fis. Mat. Natur. {\bf 34} (2023), no. 1, pp. 89-99




\end{thebibliography}
\end{document}